\documentclass[12pt,cd]{amsart}
\usepackage{amsfonts, amssymb,amsmath, amsthm, amsxtra, latexsym, mathrsfs, amscd}
\usepackage{amsmath, amssymb, amsthm, times, float, pifont, eufrak}
\usepackage{graphics}
\usepackage{pb-diagram}
\usepackage{verbatim}
\usepackage[latin1]{inputenc} %

\newtheorem{theo+}{Theorem}[section]
\newtheorem{prop+}[theo+]{Proposition}
\newtheorem{coro+}[theo+]{Corollary}
\newtheorem{lemm+} [theo+]{Lemma}
\newtheorem{deep+}  [theo+]  {Deep Result}
\newtheorem{fact+}  [theo+]  {Fact}
\theoremstyle{definition}
\newtheorem{exam+}  [theo+]  {Example}
\newtheorem{rema+}  [theo+]  {Remark}
\newtheorem{defi+}  [theo+]  {Definition}
\newtheorem{xca+}[theo+]{Exercise}

\numberwithin{equation}{section}

\usepackage{color}

\def\draft{\centerline{(Draft {\the \day}/{\the\month} \the \year.)}}
\bigskip

\def\refn#1.#2{\expandafter\def\csname#1\endcsname{[#2]}}
\def\refnr#1.{\csname#1\endcsname}

\def\Del{\Delta}

\def\fa{\mathfrak a}

\def\fg{\mathfrak g}
\def\fk{\mathfrak k}
\def\fh{\mathfrak h}
\def\fl{\mathfrak l}
\def\fm{\mathfrak m}
\def\fn{\mathfrak n}
\def\fp{\mathfrak p}
\def\fq{\mathfrak q}

\def\fsu{\mathfrak {su}}
\def\fspin{\mathfrak {spin}}
\def\fsp{\mathfrak {sp}}

\def\a{\alpha}

\def\Claminv2{|C(\Lambda)|^{-2}}

\def\Lam{\Lambda}

\def\de{d\varepsilon}

\def\Aa2D{A^{\a,2}(D)}
\def\bAa2D{\overline{A^{\a,2}(D)}}
\def\Ab2D{A^{\beta,2}(D)}
\def\bAb2D{\overline{A^{\beta,2}(D)}}

\def\Norm#1_#2{\Vert#1\Vert_{#2}}

\def\phipl12{\phi_{p_{l_1}, p_{l_2}}}
\def\phip01{\phi_{p_{0}, p_{0}}}

\def\a{\alpha}

\def\Claminv2{|C(\Lambda)|^{-2}}

\def\Del{\Delta}

\def\Lam{\Lambda}

\def\ad{\operatorname{ad}}
\def\Ad{\operatorname{Ad}}

\def\det{\operatorname{det}}

\def\exp{\operatorname{exp}}

\def\Spin{\operatorname{Spin}}
\def\SO{\operatorname{SO}}

\def\tr{\operatorname{tr}}

\def\de{d\varepsilon}

\def\Aa2D{A^{\a,2}(D)}
\def\bAa2D{\overline{A^{\a,2}(D)}}
\def\Ab2D{A^{\beta,2}(D)}
\def\bAb2D{\overline{A^{\beta,2}(D)}}

\def\phipl12{\phi_{p_{l_1}, p_{l_2}}}
\def\phip01{\phi_{p_{0}, p_{0}}}
\def\m{\underline{\bold m}}

\def\m{\underline{\bold m}}

\def\bc{\mathbb C}
\def\br{\mathbb R}

\def\bh{\mathbb H}

\def\alg/{algebra}

\def\Alg/{Algebra} 

\def\alt/{alternative} 
\def\anal/{analytic}
\def\analfunc/{\anal/\ \func/}
\def\Ans/{\it Answer. \normal}
\def\ass/{associative}
\def\nass/{non-\ass/}
\def\autom/{automorphism}
\def\homom/{homomorphism}
\def\isom/{isomorphism}
\def\bdd/{bounded}
\def\Bdd/{Bounded}
\def\bddsymdom/{bounded \sym/ \dom/}
\def\Cartdom/{Cartan \dom/}
\def\bdry/{boundary}
\def\bsd/{\bdd/ \symdom/}
\def\bv/{boundary value}
\def\cf/{{\it cf}\.}
\def\Cf/{{\it Cf}\.}
\def\charr/{character}
\def\coeff/{coefficient}
\def\comm/{commutative}
\def\cpct/{compact}
\def\compl/{complex}
\def\comp/{complex}
\def\Comp/{Complex}
\def\conf/{conformal}
\def\conj/{conjugate}
\def\conn/{connect}
\def\cont/{continuous}
\def\conv/{converge} 
\def\convc/{convergence}
\def\convt/{convergent}
\def\convx/{convex}
\def\coord/{coordinate}
\def\lcoord/{local coordinate}
\def\Corr/{Corresponding}
\def\corr/{corresponding}
\def\corrd/{correspond}
\def\cov/{covariant}
\def\decomp/{decomposition}
\def\deco/{decompose}
\def\diff/{different} 
\def\Diff/{Different} 
\def\dimn/{dimension} 
\def\distr/{distribution} 
\def\div/{diverge} 
\def\dom/{domain}
\def\eg/{\hbox{\it e.g}\.}
\def\eigenf/{eigen\-\func/}
\def\eigensp/{eigen\-space}
\def\eigenv/{eigen\-value}
\def\eq/{equation}
\def\equa/{equation}
\def\de/{\diff/ial \equa/}
\def\do/{\diff/ial operator}
\def\ode/{ordinary \de/}
\def\pde/{partial \de/}
\def\pdo/{partial \diff/ial operator}
\def\psdo/{pseudo \diff/ial operator}
\def\fin/{finite}
\def\Ex/{\it Example.\ \normal}
\def\Exnr#1/{\it Example #1.\ \normal}
\def\foll/{follow}
\def\follg/{following}
\def\Follg/{Following}
\def\func/{function}
\def\Func/{Function}
\def\Fonc/{Fonc\-tion}
\def\fonc/{fonc\-tion}
\def\Funk/{Funk\-tion}
\def\funk/{Funk\-tion}
\def\gen/{general}
\def\har/{harmonic}
\def\Hint/{\it Hint. \normal}
\def\hist/{historic}
\def\histcl/{historical}
\def\hol/{holo\-morphic}
\def\homog/{ho\-mo\-ge\-ne\-ous}
\def\hyp/{hyper\-bolic}
\def\hyperg/{hyper\-geometric}
\def\ie/{\hbox{\it i.e.}}
\def\iff/{if and only if}
\def\ineq/{inequality}
\def\infra/{{\it inf\-ra}}
\def\ultra/{{\it ult\-ra}}
\def\Inpart/{In particular}
\def\inpart/{in particular}
\def\instof/{instead of}
\def\interps/{interpolation space}
\def\interp/{interpolation}
\def\Interp/{Interpolation}
\def\interpr/{Interpretation}
\def\Intr/{Introduction}
\def\intv/{interval}
\def\inv/{invariant}
\def\invc/{invariance}
\def\Iowords/{In other words}
\def\iowords/{in other words}
\def\ipr/{inner product}
\def\irred/{irreducible}
\def\lb/{line bundle}
\def\lin/{linear}
\def\lhs/{left hand side}
\def\rhs/{right hand side}
\def\loc/{local}
\def\math/{mathematic} 
\def\mathcn/{\math/ian}
\def\manif/{manifold}
\def\meas/{measure}
\def\measl/{measurable}
\def\mero/{mero\-morphic}
\def\mon/{monomial}
\def\monog/{monogenic}
\def\mult/{multiple}
\def\multy/{multiply}
\def\multn/{multiplication}
\def\nas/{necessary and sufficient}
\def\nbd/{neighborhood}
\def\neg/{negative}
\def\nondeg/{nondegenerate}
\def\Oohand/{On the other hand}
\def\oohand/{on the other hand}
\def\Oonhand/{On the one hand}
\def\oonhand/{on the one hand}
\def\oper/{operator}
\def\orth/{ortho\-gonal}
\def\orthon/{ortho\-normal}
\def\otoh/{on the other hand}
\def\quat/{quaternion}
\def\pp/{\hbox{a. e.}}
\def\psh/{plurisubharmonic}
\def\pol/{polynomial}
\def\pot/{potential}
\def\pos/{positive}
\def\princ/{principle}
\def\prob/{probability}
\def\proj/{projective}
\def\projn/{projection}
\def\Proof/{\it Proof:\normal}
\def\Rem/{\it Remark\normal}
\def\Remnr#1/{\it Remark\ \normal #1. }
\def\rep/{representation}
\def\reps/{representations}
\def\meta/{metaplectic representation}
\def\repr/{reproducing}
\def\reprker/{reproducing kernel}
\def\resp/{respective} 
\def\resply/{respectively}
\def\restr/{restriction}
\def\sa/{self-adjoint}
\def\st/{such that}

\def\sol/{solution}
\def\ru/{space}
\def\sph/{spherical}
\def\ssp/{sub\ru/}
\def\sym/{symmetric}
\def\Sym/{Symmetric}
\def\symb/{symbol}
\def\symbc/{symbolic}
\def\symdom/{\sym/ domain}
\def\symp/{symplectic}
\def\Theor#1/{\fet Theorem #1.\ \normal}
\def\Lem#1/{\fet Lemma #1.\ \normal}
\def\Lemma/{\fet Lemma.\ \normal}
\def\topl/{topology}
\def\topll/{topological}
\def\transf/{transform}
\def\transl/{translation}
\def\transfn/{transformation}
\def\transv/{transvectant}
\def\trig/{trigonometric}
\def\tril/{trilinear}
\def\trilf/{trilinear form}
\def\uhp/{upper halfplane}
\def\uhs/{upper halfspace}
\def\vb/{vector bundle}
\def\vf/{vector field}
\def\vsp/{vector space}
\def\wrt/{with respect to}
\def\Wlog/{Without loss of generality}
\def\a{\alpha}

\def\Lam{\Lambda}

\def\Ab/{Abel}
\def\Ban/{Banach}
\def\Bansp/{\Ban/ space}
\def\Belt/{Bel\-tra\-mi}
\def\Berg/{Berg\-man}
\def\Bern/{Ber\-nou\-lli}
\def\Berz/{Berezin}
\def\Bess/{Bessel}
\def\Cart/{Car\-tan}
\def\Cay/{Cay\-ley}
\def\CG/{Clebsch-Gordan}
\def\Cl/{Clifford}
\def\CR/{Cauchy-Rie\-mann}
\def\Dir/{Dirichlet}
\def\Eucl/{Euclide}
\def\Eucln/{Euclidean}
\def\F/{Fourier}
\def\Hank/{Hankel}
\def\Hankf/{\Hank/ form}
\def\Herm/{Hermite}
\def\Hilb/{Hilbert}
\def\Hilbs/{Hilbert space}
\def\Hilbsp/{Hilbert space}
\def\HS/{Hilbert-Schmidt}
\def\Lag/{La\-grange}
\def\Lap/{La\-place}
\def\LapBelt/{\Lap/-\Belt/}
\def\Leb/{Lebesgue}
\def\Marc/{Mar\-cin\-kie\-wicz}
\def\Moeb/{Moebius}
\def\Moebt/{Moebius transformation}
\def\Moebtransfn/{Moebius transformation}
\def\Pla/{Plan\-che\-rel}
\def\Poin/{Poin\-car\'e}
\def\Riem/{Rie\-mann}
\def\Riemn/{\Riem/ian}
\def\psRiemn/{pseudo-\Riem/ian}
\def\Riems/{Rie\-mann surface}
\def\Schroe/{Schr\"odinger}
\def\Weier/{Weier\-strass}
\def\anal/{analytic}
\def\bsd/{bounded symmetric domain  }
\def\bdd/{bounded}
\def\calc/{calculation}\def\conj{conjugate}
\def\calci/{calculating}\def\eg{e.g.}
\def\conj/{conjugate}
\def\deco/{decomposition}
\def\eg/{e.g.}
\def\fct/{function}
\def\gp/{group}
\def\hw/{highest weight}
\def\hwv/{highest weight vector}
\def\hwvs/{highest weight vectors}
\def\lw/{lowest weight}
\def\lwv/{lowest weight vector}
\def\lwvs/{lowest weight vectors}
\def\hds/{holomorphic discrete series}
\def\iff/{if and only if}
\def\inv/{invariant}
\def\irrde/{irreducible decomposition}
\def\meas/{measure}
\def\transf/{transform}
\def\rep/{representation}
\def\resp/{respectively}
\def\inters/{intertwines}
\def\interg/{intertwining}
\def\meta/{metaplectic representation}
\def\qu/{quaternion}
\def\rep/{representation}
\def\symdom/{ symmetric domain}
\def\st/{such that}
\def\shd/{subhead}
\def\transf/{transform}
\def\wrt/{with respect to}

\def\Norm#1#2#3{\Vert#1\Vert^{#3}_{{#2}}}

\def\tr{\operatorname{tr}}

\begin{document}
\title[Restriction of unitary representations
]{
Restriction  to symmetric subgroups
of unitary representations of rank one 
semisimple Lie groups
 }
\author{B. Speh and G. Zhang}
\address{B. Speh, Department of Mathematics, Cornell University,
  NY14853-4201, USA.  
\newline
Email: 
speh@math.cornell.edu}

\address{G.Zhang, Mathematical Sciences, Chalmers University of Technology and
Mathematical Sciences, G\"oteborg University, SE-412 96 G\"oteborg,
Sweden. 
\newline
Email:  genkai@chalmers.se}

\thanks{Research by B.Speh partially supported by NSF grant DMS-0901024  and research by G. Zhang partially supported by 
the Swedish
Science Council (VR)}
\begin{abstract}
We consider the spherical complementary
series of rank one Lie groups
$H_n=\SO_0(n, 1; \mathbb F)$
for $\mathbb F=\mathbb R, \mathbb C, \mathbb H$. We prove that there
exist finitely many discrete components
in its restriction under the subgroup
$H_{n-1}=\SO_0(n-1, 1; \mathbb F)$.
This is proved by imbedding
the complementary series 
into analytic continuation
of holomorphic 
discrete series
of $G_n=SU(n, 1)$, $SU(n, 1)\times SU(n, 1)$
and $SU(2n, 2)$ and by the branching
of holomorphic representations
under the corresponding subgroup
 $G_{n-1}$.
\end{abstract}

\maketitle

\baselineskip 1.35pc

\section{Introduction}

This paper is a continuation
of  \cite{Speh-Venk-2}
by Speh-Venkataramana   and 
\cite{gz-brch-rk1} by Zhang.  Consider the complementary series representations  $(\pi_{\mu}, H)$ of the 
 Lie group $H=SO_0(n, 1; \mathbb F)$
for $\mathbb F=\mathbb R, \mathbb C, \mathbb H$ and its restriction to the subgroup  $H_{1}=SO_0(n-1, 1; \mathbb F)$.
It is proved there that for certain parameters $\mu$ the restricted representation has a direct summand which is isomorphic to a complementary series representation of the smaller group $H_1$;
similar results are proved
also for the exceptional Lie group
$H=SO_0(2, 1; \mathbb O):=F_{4(-20)}$,
and $H_1=Spin(8, 1)$. 

In  \cite{Speh-Venk-2} the complementary series  $(\pi_{\mu}, H)$ 
of
$H=SO_0(n, 1; \mathbb R)$
are realized in its non-compact picture,
namely as appropriate spaces
of distributions on the space $\mathbb R^{n-1}$
and it is proved that, roughly speaking, 
up to the Fourier transform,
the trivial extension of distributions
on the subspace $\mathbb R^{n-2}$ 
to distributions $\mathbb R^{n-1}$
defines an $H_1$--intertwining
operator and thus determines for proper values of $\mu$
a discrete summand of $(\pi_{\mu}, H)$. 

The general rank one case is treated in 
\cite{gz-brch-rk1} 
 using 
the compact picture realizing
the complementary series as certain Sobolev
type spaces on the sphere in $\mathbb F^n$, 
the discrete components being 
 obtained by restricting the functions
to the sub-sphere in  $\mathbb F^{n-1}$.
This requires some detailed estimates
for the branching of representations
of the maximal compact subgroup $K\subset H$
 on the space of spherical
harmonics on $\mathbb F^{n}$
under the  compact subgroup $L\subset H_1$.

\medskip
In the present paper we determine additional discrete summands of the restriction of  $(\pi_{\mu}, H)$ to $H_!$. These are also isomorphic to complementary series representations of $H_1$ and we determine explicitly  the parameters of these complementary series representations of $H_1$. 
See theorem 5.1 for a precise statement.

\medskip
In the proof we shall use a rather different idea than in the previous papers. 
We consider $H$ respectively $H_1$ as subgroups of a larger groups $G$, respectively $G_1$ so that we have the following
commutative diagram of subgroup inclusions
\[
\begin{diagram}
\node{G} \node{G_{1}}\arrow{w}
\\
\node{H} \arrow{n} \node{H_{1}}\arrow{n}\arrow{w}
\end{diagram}
\]
Here we shall 
consider the Hermitian groups $G=SU(n, 1)$,  $
SU(n, 1)\times
SU(n, 1)$,  $SU(2n, 2)$.




We start with an 
scalar holomorphic representation $(U_\nu, G)$
of $G$ realized on a Hilbert space
of holomorphic functions on the unit ball $G/K$,
and consider its branching
along the diagram. Denote $(U_{\nu}^\flat, G_{1})$
the scalar holomorphic representations of $G_{1}$. 
The restriction of $(U_\nu, G)$ to $G_1$ is a direct sum of unitary representations  and the branching rules are obtained 
by  normal holomorphic  differentiation
 $\mathcal D^k$ of degree $k$
along  the submanifold manifold $G_{1}/K_{1}$ in  $G/K$. 
See also  \cite{kobayashi-06}
for an abstract study. The representations
appearing in the restriction will be also holomorphic representations
and we are interested in the scalar ones
 $(U_{\nu'}, G_{1})$.
The branching
of $(U_\nu, G)$ under $H$ involves
more analytic issues, and has been studied
extensively   \cite{Neretin-plan-beredef,  Dijk-Hille-jfa, 
gz-bere-rbsd, gz-br2, gz-adv}. In particular
it is proved that
the complementary series
 $( \pi_{\mu}, H)$ 
appears in  $(U_\nu, G)$ for certain pairs $(\mu, \nu)$.
In the compact-picture of complementary
series as distributional spaces
on the boundary $S=H/MAN$ of the unit ball $H/L$
there is an  intertwining operator
$$
J_\nu:
(\pi_{\mu}, H)\to {(U_\nu, G)}
$$
obtained by pairing the distributions
on $S$
in the complementary series with the reproducing
kernel of the holomorphic representation
as a function on $S \times \, G/K$, resulting in 
a kind of  holomorphic extension of complementary
series $(\pi_{\mu}, H)$  in 
$(U_\nu, G)$.
The same
is true also for 
 $(\pi_{\mu'}, H_{1})$ 
and  $(U_{\nu'}, G_{1})$. The problem reduces thus to
proving  the operator
$(J_{\nu'}^\flat)^\ast
\mathcal D^k
J_{\nu} $ in the diagram
\[
\begin{diagram}
\node{(U_\nu, G)}
 \arrow{e,t}{\mathcal D^k} 
\node {(U_{\nu'}, G_1)}
\arrow{s,r}
{
(J_{\nu'})^\ast
} 
\\
\node{(\pi_{\mu}, H)}
\arrow{e}
\arrow{n,l}
{
J_{\nu}
} 
\node
{(\pi_{\mu'}, H_1)}
\end{diagram}
\]
is non-zero for appropriate $\mu$ and $\mu'$. 
This will be done by computing the operator acting on certain vector 
in $(\pi_{\mu}, H) $; see Theorem 5.1 below.

It is clear 
that,
compared with the method in \cite{gz-brch-rk1} 
where we  study boundedness property of the restriction
by $K$-type expansion,  there
are   some advantages
with  the present approach as all the intertwining 
operators involved have been studied before and
there also is flexibility of choosing the degree
of differentiation. 


The above procedure of holomorphic extension
and differentiation can be applied 
to general spherical (not necessarily unitary)
representations and we can 
find formally intertwining operators from 
spherical principal series
of $H=SO_o(n, 1; F)$ to that of 
$H_1=SO_o(n-1, 1; F)
$. In the above compact realization 
of the representations
as distribution spaces on the 
sphere $S$ for $H$ and respectively
on the subsphere $S^\flat$  for $H_1$, 
those operators will be 
polynomials of normal differentiations
orthogonal to the totally geodesic submanifold
$S^\flat \subset S$ and of horizontal
differentiations parallel to 
$S^\flat \subset S$.
The exact formulas can be rather involved.
For the real case,
$H=SO_o(n, 1)$, $H_1=SO_o(n-1, 1)$, those operators
have been 
found by Juhl \cite{Juhl-book2}
using  computations
involving enveloping algebras
and Verma modules. See also Kobayashi \cite{kobayashi-helgason85}.
In the last section we shall give for $SO_0(n,1)$
a direct proof of our result
using these differential intertwining operators in 
the non-compact realizations.

The results in  \cite{Speh-Venk-2} suggest that similar results should be true for the nonspherical complementary series representations of $SO_o(n,1)$ and more generally for complementary series of all connected semi simple groups of real rank 1.

\medskip
The branching of complementary series is of interests
  in
  automorphic
representations \cite{Speh-Venk-2},  \cite{Orsted-Speh}, 
\cite{Burger-Li-Sarnak}, \cite{Burger-Sarnak},
\cite{Bergeron-imrn}. 

After a preliminary version of this paper was completed
we were informed of the preprint \cite{Mollers-Oshima}
where they obtained a complete irreducible decomposition  of
$(\pi_{\mu}, SO(n, 1))$ under the subgroup $SO(n-1, 1)$. The case $n=3$
has  been obtained earlier in mathematical physics litteratures; see
e.g. \cite{Mukunda}.

\medskip
Part this work
was started during an AIM work ``Branching
problems in unitary representations'',
MPIM, Bonn, July 2011. We thank AIM
and MPIM
for providing a very stimulating
environment. We are also grateful
to Roger Howe and T. N. 
Venkataramana for several
conversations,  and  Bent \O{}rsted for drawing our attention
to the work \cite{Juhl-book2}.

\section{Spherical representations of  connected semi simple Lie groups of real rank 1}

We recall very briefly some known results and fix
some notation.
Let $\mathbb F=\mathbb R, \mathbb C, \mathbb H$
be the real, complex and quaternionic numbers
of real dimension $\iota=1, 2, 4$.
Let $\mathbb F^{p, q}:=\mathbb F^{p+q}$
be equipped with 
the quadratic form $|x_1|^2 +\cdots +|x_p|^2 -|x_{p+1}|^2
-\cdots |x_{p+q}|^2$
for signature $(\iota p, \iota q)$.
Denote $SO_0(p, q; \mathbb F)
$ the connected component
of group of $\mathbb F$-linear transformations
on $\mathbb F^{p, q}$ preserving
the  form and having determinant $1$, 
with $\mathbb F$ acting on the right. In the usual notation it is equal to $SO_0(p, q),
SU(p, q), Sp(p, q)$ respectively.
Elements in $SO_0(p, q; \mathbb F)$ and 
its Lie algebra 
$\mathfrak{so}(p, q;
\mathbb F)$
  will be written
as $(p+q)\times (p+q)$ block $\mathbb F$-matrices 
$$
\begin{bmatrix} a_1&a_2\\
a_3&a_4
\end{bmatrix}
$$
where $a_1, a_2, a_3, a_4$ are of size $p\times p,
p\times q, q\times p, q\times q$ respectively.
We shall consider  subgroups
 $SO_0(r, s; \mathbb F)$
in 
$SO_0(p, q; \mathbb F)$, 
 $r\le p, s\le q$,
and fix the realizations as follows.
Write $\mathbb F^{p, q}=\mathbb F^{r} \oplus
\mathbb F^{p-r} \oplus \mathbb F^{p} \oplus
\mathbb F^{q-s}$.
 The group
 $SO_0(r, s; \mathbb F)$
 consists of elements
in $SO_0(p, q; \mathbb F)$ fixing
all vectors $\mathbb F^{p-r} \oplus
\mathbb F^{q-s}$ in the above decomposition.
The 
Lie algebra 
$\mathfrak{so}(r, s;
\mathbb F)$ consists of all elements in 
$\mathfrak{so}(p, q;
\mathbb F)$ with the entries $a_j$ being diagonal matrix
$$
a_j=\begin{bmatrix} b_j&0\\
0&0
\end{bmatrix}, 
$$
with the obvious sizes of the blocks.
The Riemannian symmetric space 
of $SO_0(p, q; \mathbb F)$ will
be 
realized
as the matrix ball in $M_{p, q}(\mathbb F)$,
$$
SO_0(p, q; \mathbb F)/S(SO(p; \mathbb F)\times
SO(q;\mathbb F))=\{z\in M_{p, q}(\mathbb F); zz^\ast < I\}.
$$ 
The above identification defines
a  totally geodesic
submanifold
$$
SO_0(r, s; \mathbb F)/
S(SO(r; \mathbb F)\times
SO(s;\mathbb F))
\subset  
SO_0(p, q; \mathbb F)/
S(SO(p; \mathbb F)\times
SO(q;\mathbb F)).
$$ 

\medskip

We  fix now
$H=SO_0(n, 1; \mathbb F)$.
The maximal compact subgroup $L$ is equal to 
$SO(n), \linebreak S(U(n)\times U(1)),
Sp(n)\times Sp(1)$ respectively and
and $H/L$
is a Riemannian symmetric space of rank one.

Let $\fh=\fl +\fq$ be the corresponding
Cartan decomposition of the Lie algebra $\fh$.  Elements  in $\fq$
are of the form
$$
\begin{bmatrix} 0&v\\
v^\ast&0
\end{bmatrix} 
$$
with $v\in \mathbb F^n$. Here  $v^\ast=
\overline{v^T}$ the conjugate transpose where
$\overline{a}$ is the conjugation in $\mathbb F$. We
shall thus identify $\fq$ with $\mathbb F^n$.
We fix
$$
H_0=e_1\in \fq=\mathbb F^n, e_1^T=[1, \cdots, 0]
$$
for $\mathbb F=\mathbb R$ or $ \mathbb C$;
for $\mathbb F=\mathbb H$ we let $H_0$
be as in \cite{gz-brch-rk1} 
so that $\Ad(H_0)$ has eigenvalues $\pm 2,
\pm 1, 0$.
Then $\fa:=\mathbb RH_0$
is a maximal abelian subspace of 
$\fq$. The root space decomposition of $\fh$ under  $H_0$ 
is 
$$
\fh=\fn_{-2} + \fn_{-1} +
(\fa+ \fm) + \fn_{1} +\fn_{2} 
$$
with $\pm 2, \pm 1, 0$,
if $\mathbb F=\mathbb C, \bh$,
and with the convention that $\fn_{2} =0$
if $\mathbb F=\mathbb R$.
  Here $
\fm\subset \fl$ is the zero root space.
We denote by $\fn=\fn_1 \oplus \fn_2$ the sum of the positive root spaces.
 Then
$\fm+\fa +\fn$
is a maximal parabolic subalgebra
of $\fh$.

Denote $M, A, N$ 
the corresponding 
subgroups with Lie algebras
$\fm, \fa, \fn$. Then  $M=SO(n-1), SU(n-1),$  respectively $Sp(n-1)\times Sp(1)$
and $MAN$ is a maximal parabolic
subgroup of $H$.  We normalize the $K$-invariant norm on $\fp$
so that $H_0$ is a unit vector.
The space $S=K/M$ is then the unit sphere
in $\fp$ with $M$ the isotropic subgroup of $H_0\in S$.

Let $\rho=\rho_{H}
$
be the half sum of positive roots. Then
$$
\rho(H_0)=
\begin{cases}
\frac{n-1}2, &\quad \mathbb F=\br\\
n, &\quad \mathbb F=\bc\\
2n+1, &\quad \mathbb F=\bh
\end{cases}
$$
and we shall identify $\rho=\rho(H_0)$.
We let $H_{1}\subset H$ be the
 the subgroup  $SO_0(n-1, 1, \mathbb F)$ with the above
convention of the notation.
Let $\fh_1=\fl_1+\fq_1$
be the Cartan decomposition of $\fh_1$
obtained
under the Cartan involution on $\fh$
restricted to $\fh_1$. 
The subspace  $\fa=\mathbb RH_0\subset \fq$
is also in $\fq_1$,  and
we have the corresponding root space decomposition
$$
\fh_1=\fn_{-2, 1} + \fn_{-1, 1} +
(\fa+ \fm\cap \fl_1) + \fn_{1, 1} +\fn_{2, 1}, $$
with
$$ \fn_{\pm 2, 1}=\fn_{\pm 2},\ \  \
\ 
\fn_{\pm 1, 1}=\mathbb F^{n-2} \subset 
\fn_{\pm 1}=\mathbb  F^{n-1}, 
$$
and $\fn_1/\fn_{1, 1}=\mathbb F$.

The corresponding subgroup $P_1= P\cap H_1=(L_1\cap M)AN_1$
is a maximal parabolic subgroup of $H_1$,
where $N_1$ is the nilpotent
subgroup with Lie algebra $\fn_1=\fn_{1, 1} +\fn_{2, 1} $.
The homogeneous
space 
$L/K\cap M$
is a  subsphere  
$S_1:=S^{\iota (n-1)-1}$ 
in the sphere
$S:=K/M
= S^{\iota n-1}$,
\begin{equation}
  \label{eq:1-S-var}
  S_1=L/K\cap M\subset S=K/M.
\end{equation}
 Here $\iota = 1,2$ or 4.

\medskip

For $\mu\in \mathbb C$
let $\pi_{\mu}$ 
be the spherical principal series representation of $H$
induced  from $MAN$ consisting of measurable functions $f$ on $H$
such that
\begin{equation}
  \label{eq:ind-norm}
f(g me^{tH_0}
n
)=
e^{-\mu t}f(g), \ \ m\,  e^{tH_0}n\in MAN  
\end{equation}
and $f\big{|}_{L}\in {\it L}^2(L)$. 
In particular elements $f\big{|}_{L}$
are invariant under the right action of  $M$
and thus $f\big{|}_{L}\in {\it L}^2(L/M)={\it L}^2(S)$. 
The group action of $\pi_{\mu}(h)$ of $H$ on ${\it L}^2(S)$ is
\begin{equation}
  \label{eq:H-act}
\pi_{\mu}(h) f(u)=|cu +d|^{-\mu} f((au+b)(cu+d)^{-1}), \quad
h^{-1}=
\begin{bmatrix} a&b\\
c&d
\end{bmatrix}
\end{equation}
for $f\in L^2(S)$.
This can be derived for example by using that if
\[h=
\begin{bmatrix} a&b\\
c&d
\end{bmatrix}\in H, \, u=k\cdot H_0=kMAN\in S=L/M.\]
then
\begin{equation}
  \label{eq:H-jacobi}
e^{-2\rho A(h k)}
=J_{h}(u) =|cu +d|^{-2\rho},  \quad 
\end{equation}
as proved in \cite{Johnson-Wallach, CDKR-jga}. 
Here $e^{H(h)}$ is the $A$-component
in the Iwasawa decomposition $g=ke^{A(h)}n$
of $H=LAN$ of $H$, $J_g(u)$ is the Jacobian of $g: u\to gu$
on $S$.

The Lie algebra action $ (\pi_{\mu}, \fh)$
on the space ${\it L}^2(S)_L$
of $L$-finite elements in ${\it L}^2(S)$
defines a Harish-Chandra module which we denote also denote $\pi_{\mu} $.

 If  $\mu $ is real and $0<\mu<\rho(H_0)$ provided that $\mathbb F=\mathbb R, \mathbb C$
 and $2 <\mu <\rho(H_0)$  for $\mathbb F= \mathbb H$
 the representations $({\it L} ^2(S)_L,  \pi_{\mu}, \fh)$ are unitarizable and 
 called complementary series.
The unitarization of $({\it L}^2(S)_L,  \pi_{\mu}, \fh)$
can be done through the Knapp-Stein intertwining operator.
Its exact formula  will not be
used 
in our paper, but we recall it briefly
as it motivates our definition of the $J$-operator
in Theorem 4.2.

Let $\langle x, u\rangle
=x_1\bar u_1 +\cdots +x_n\bar u_n$ be the $\mathbb F$-valued
inner product on $\mathbb F^n$.
Let 
\begin{equation}
\label{eq:Lam-J}
\Lambda_{\mu}f(x)=
C_{\mu}\int_S \frac{f(u)}
{|1-\langle x, u\rangle|^{\mu}}
 du.
\end{equation}
We note that the integral is absolutely convergent for $\mu<\rho$
for any bounded function $f$ on $S$.
As an operator on $L$-finite elements in ${\it L}^2(S)$ it
has meromorphic continuation in $\mu$ to the complex plane.
Then
\begin{equation}
  \label{eq:kn-st}
\Lambda_{2\rho-\mu}: ({\it L}^2(S)_L, \pi_{\mu})
 \to  ({\it L}^2(S)_L, \pi_{2\rho-\mu})
\end{equation}
is a $\fh$-intertwining operator.
A $\fh$-invariant Hermitian form on ${\it L}^2(S)_L$ is
given by
$$
(f, g)_\mu=(\Lam_{2\rho-\mu}f, g).
$$
The decomposition of ${\it L}^2(S)$ under $L$
is multiplicity free  \cite{Johnson-Wallach},
$L^2(S)=\sum_{\tau}W^\tau$.  In particular
the operator $\Lam_{2\rho-\mu}$
is diagonal under the decomposition,
$$
\Lam_{2\rho-\mu}
 f =\lambda_{\mu}(\tau) f, \quad  \mbox{ for }  f\in W^\tau,
$$
with ${\lambda_\mu(\tau)}$ explicitly
found in \cite{Johnson-Wallach}. 
\begin{theo+} \cite{Kostant-BAMS} \ \ 
 If
 \begin{enumerate}
\item 
$\mathbb F=\br$, $0 <\mu < n-1=2\rho_{H}
$;
\item  $\mathbb F=\bc$,  $0 <\mu < 2n=2\rho_{H}$;
\item  $\mathbb F=\bh$, $2 <\mu < 4n
=2\rho_{H}-2$,
\end{enumerate}
there is a positive definite
$\fh$-invariant form on $L^2(S)_L$
defined on the isotopic component $W^\tau $ by
$$
\Vert w\Vert_{\mu}^2={\lambda_\mu(\tau)}
\Vert w\Vert^2.
$$
Its completion  $(\mathcal C_\mu, \pi_{\mu}, H)$
defines an unitary irreducible complementary series representations
of $H$.\end{theo+}

The unitary
representations $\pi_{\mu}$ and $\pi_{2\rho-\mu}$
are unitarily equivalent. To avoid possible
confusions we always assume that $\mu<\rho_{H}$, i.e.
$\mu$ is in
the left half of the above interval for
complementary series.

\medskip

\section{Analytic continuation 
 of 
holomorphic discrete series
of  unitary groups and their restrictions 
to a unitary subgroup}

We fix  pairs $(G, G_1)$ as follows
\begin{enumerate}
\item
$\mathbb F=\mathbb R$, $(G, G_1)=(SU(n, 1), SU(n-1))$;
\item $\mathbb F=\mathbb C,$ 
$(G, G_1)=(SU(n, 1)\times
SU(n, 1), SU(n-1, 1)\times
SU(n-1, 1))$;
\item
$\mathbb F=\mathbb H,$
 $(G, G_1)=(SU(2n, 2), SU(2(n-1), 2))$.
 \end{enumerate}
 In this action we discuss  the branching
of analytic continuation
of the holomorphic
discrete series
of $G$  under under restriction to $G_1$. Results for  $\mathbb F=\bc$,
i.e., the tensor product case,  are trivial generalization
of those for  $\mathbb R$ and will be treated
only very briefly.

\medskip

The group $G$  is a Hermitian Lie group
and its symmetric space $G/K$ will be realized
as a bounded symmetric domain $D$.

 More precisely
we view the complex matrices $z$ in 
 $$M_{2n, 2}(\bc)=M_{2, 2}(\bc)
\oplus \cdots \oplus M_{2, 2}(\bc)
$$ as  column vectors $z=[z_1, \cdots, z_n]^T$ with 
each $z_j$ being $2\times 2$-matrices and define
$$
V=\begin{cases} \bc^n & \text{Case} \, \mathbb R
\\
\bc^n \oplus \overline{\bc^n} & \text{Case} \, \mathbb C
\\
 M_{2n, 2}(\mathbb C) & \text{Case}\, \mathbb H.
\end{cases}
$$
The group $G$ acts on $$D=\{z\in V; zz^\ast <I\}$$ by
$$
gz=(az+b)(cz+d)^{-1}, \,  
g=
\begin{bmatrix} a&b\\
c&d
\end{bmatrix}\in G,
$$
in the cases of $\mathbb R$  and $\mathbb H$,
it acts on the second component $(z_1, z_2)\in D$ 
as anti-holomorphic mappings in the case $\mathbb C$.

\medskip

The holomorphic discrete series
of $G$ are realized
as weighted Bergman spaces of
holomorphic functions on $D$ and are 
determined by their reproducing kernels. Let
$$
h(z, w)=
\begin{cases} 1-\langle z, w\rangle & \text{Case }\, \mathbb R
\\
(1-\langle z_1, w_1\rangle )
(1-\langle w_2, z_2\rangle )
& \text{Case }\, \mathbb C
\\
\det(1- w^\ast z) & \text{Case }\,\mathbb H.
\end{cases}
$$
Note that  the $h$-function for  ${\mathbb F}=\br$ can also be written
as $\det(1-w^\ast z)$.

When $\nu$ is in the set
$$
\begin{cases} (n, \infty) & \text{Case }\, \mathbb R, 
 \mathbb C \\
(2n+1, \infty) & \text{Case}\, \mathbb H
\end{cases}
$$
the kernel $h(z,w)^{-\nu}$
is the reproducing kernel of a weighed Bergman
space $\mathcal H_\nu$ 
with a probability measure
on $D$, i.e. the define discrete series representations of the universal covering
of $G$.

 The group $G$
acts on  $\mathcal H_\nu$ 
via
$$
U_\nu (g) f(z)= \det(cz+d)^{-\nu} f(g^{-1}z), \quad
g^{-1}=
\begin{bmatrix} a&b\\
c&d
\end{bmatrix},
$$
defining a (projective) unitary representation
$(U_\nu, G)$; in the case of $\mathbb C$ it is
the tensor product.
 The set of all $\nu\ge 0$
for which 
 $h(z,w)^{-\nu}
$ is positive 
definite (also called Wallach set) is well-known for general bounded
symmetric domains \cite{FK-book, FK, Rossi-Vergne-acta,
  Wallach-tams}; it consists of an open interval
and  discrete points and we shall only need open interval.
In our parametrization  we have

 \begin{lemm+} The functions $h(z, w)^{-\nu}$
is positive definite if
$\nu $ is in the set
$$
\begin{cases} (0, \infty) & \text{Case }\, \mathbb R, \, 
\mathbb C, \, 
\\
(1, \infty) & \text{Case }\,  \mathbb H.\\
\end{cases}
$$
and it determines  a unitary  (projective)
representation $(U_\nu, G)$
on a Hilbert space of
holomorphic functions on $D$ with reproducing kernel
 $h(z, w)^{-\nu}$.
In particular the representation  $(U_\nu, G)$
is unitary and not  in the discrete series if
$\nu$ is in the set
$$
\begin{cases} (0, n) & \text{Case }\, \mathbb R, \, 
\mathbb C, \, 
\\
(1, 2n+1) & \text{Case }\,  \mathbb H\\
\end{cases}.
$$
\end{lemm+}

\medskip
We denote the analytic continuation also
by $
(\mathcal H_\nu, U_\nu, G)=
(\mathcal H_\nu(D), U_\nu, G)
$.
In this paper we are concerned with representations $
(\mathcal H_\nu, U_\nu, G)$ whose parameter
 $\nu$ is in a subinterval of  the above interval, i.e.,it is not a 
discrete holomorphic representations.

 For the case $\mathbb C$ with $G=SU(n, 1)\times SU(n, 1)$
we let
$\mathcal H_\nu(D\otimes \bar D)
=\mathcal H_\nu(D)\otimes \overline{\mathcal H_\nu( D)}$
and $(\mathcal H_\nu(D\otimes \bar D),
U_\nu, G)$ be the tensor product of $(U_\nu\otimes\overline{U_\nu}, 
SU(n, 1)\times SU(n, 1))$, for $\nu>0$.

\medskip

To find the branching of $(\mathcal H_\nu, U_\nu, G)$
under $G_1$ we realize $D_1=G_1/K_1$ as a submanifold
of $D=G/K$ as in \S2.
For any polynomial $p$
we denote $\partial_p$ the
corresponding differential operator.
Let 
$$
Q=
\begin{cases} 
z_n,  & \text{Case} \, \mathbb R
 \\
z_n \bar {w_n},
&\text{Case} \, \mathbb C \\
\det(z_n),  & \text{Case}\,
 \mathbb H.
\end{cases}
$$
The polynomial $Q$ on $V$, as tangent space of $D=G/K$, is
invariant under the subgroup $K_1\subset K$.

\begin{theo+}  Suppose that $\nu$ satisfies the assumptions of 3.1.
\begin{enumerate}
\item The restriction of  of $(U_{\nu}, G)$
under $G_1$ is a discrete sum of  holomorphic
representations of $G_1$.
\item The holomorphic representation $(U_{\nu +j}, G_{1})$ for a nonnegative integer j,
is a discrete summand of  representation $(U_{\nu}, G)$, restricted to $G_1$
\item A non-zero $G_1$-intertwining operator 
from  $(U_{\nu}, G)$ to $(U_{\nu+j}, G_{1})$
is given
by
$$
\mathcal D^j: f\to \partial_{Q^j}f\large{|}_{D_1}.
$$
\end{enumerate}
\end{theo+}

\noindent {\it Proof.} \ 
Since $(U_{\nu}, G)$ is a highest weight module the first assertion follows from 
\cite{kobayashi-06}.    Following  the ideas in \cite{jacobsen-vergne}  to prove the second and third assertion
we realize in the case $\mathbb R$
the group $G_1\subset G$ acts on $\mathcal H_\nu$
via
$$
U_\nu(g) f(z', z_n)= (cz'+d)^{-\nu} f(
(az'+b)(cz'+d)^{-1}, 
z_n(cz'+d)^{-1}), \quad g^{-1}=
\begin{bmatrix} a&b\\
c&d
\end{bmatrix}
\in G_1.
$$
Thus $(\partial_n^j)_{z_n=0}$ is an intertwining operator
$$
\partial_n^j
( U_\nu(g) f)(z', 0)= (cz'+d)^{-\nu-j}
 \left(\partial_n^j f\right)
((az'+b)(cz'+d)^{-1}, 0)
=\pi_{\nu+j}(g) \left (\partial_n^j f(\cdot, 0)
\right).
$$
The case $\mathbb C$
follows trivially. 
In the case $\mathbb H$ the intertwining property
is a simple consequence of
$$
\det(\partial_n)^j
 (f(Az_n))
=\det(A)^j (\det(\partial_n)^j f)(Az_n)
$$
for any function $f$ on $M_{2, 2}(\mathbb C)$.
$\Box $


\section
{The branching of  $(U_\nu, G)$
under  $H$. The appearance of the complementary series 
$(\pi_{\mu}, H)$}

 We consider $H\subset G$, $H_1\subset G_1$,
written as $(H, H_1)\subset (G, G_1)$,
where \linebreak
$(H, H_1)=(SO_0(n, 1; \mathbb F),
SO_0(n-1, 1; \mathbb F))$
\medskip
\begin{itemize}
\item If $\mathbb F = \mathbb R$ then $(H, H_1)
\subset (G, G_1)$ is the fix point set of the unitary groups under complex conjugation.\\
\item If   $\mathbb F=\mathbb C$ 
the subgroup pair $(H, H_1)\subset (G, G_1)$ 
 is the diagonally embedded. \\
\item If  $\mathbb F=\mathbb H$ 
we first 
identify the quaternions $\bh$ 
as the real form 
$$\bh=\{\begin{bmatrix}
 a+ib & c+id\\
-c+id & a-ib
\end{bmatrix}\subset M_{2, 2}(\bc)
$$
in $M_{2, 2}(\bc)$,
 i.e.  $\mathbb F=\mathbb H$ 
is the set of elements in  $M_{2, 2}(\bc)$
fixed by the real involution
 $J: z= x +iy\to x -iy$ in $M_{2, 2}(\bc)$.  Observe
that under this identification
$$
\det \bh =|\bh|^2.
$$
The involution
defines a corresponding $J$
on $M_{2+2n, 2+2n}(\bc)$. The group
$H=Sp(n, 1)$ is then the 
realized as the set of elements in
$SU(2n, 2)\subset M_{2+2n, 2+2n}(\bc)$ fixed
by $J$.
\end{itemize}

In this notation the complex-valued
polynomial $\det: M_{2,2}(\bc) \rightarrow \bc $ is the unique $Sp(1)\times Sp(1)$-invariant
holomorphic polynomial  of degree $2$.
The complex matrices $z$ in 
 $M_{2n, 2}(\bc)$ are then complexified
quaternions column vectors $z=[z_1, \cdots, z_n]^T$. Consider  the vectors
$z=[z_1, \cdots, z_n]^T$ with quaternionic entries, i.e.,
in the real space $\mathbb H^n$.
Then 
$$
p(z)=\det z_1 +\cdots +\det z_n =
| z_1|^2 +\cdots +|z_n|^2 =\Vert z\Vert^2,
$$ is invariant under 
$Sp(n)\times Sp(1)$ with 
$Sp(n)$ acting on the left and
$Sp(1)$  on the right. 
Thus 
the holomorphic polynomial $p(z)=\det z_1 +\cdots +\det z_n$ is also invariant
on the complexification $M_{2,2}(\bc)$ since
$\mathbb H^n$ is a totally real subspace in $M_{2,2}(\bc)$.

\medskip
The main result in the section is

\begin{theo+}
Every complementary series $(\pi_{\mu}, H)$,
$\mu <\rho_H$,
of $H$ is a direct summand of a holomorphic representation of $G$.
\end{theo+}


\medskip  To obtain the  precise formulas for the branching
we now recall the branching of  $(\mathcal H_\nu, U_\nu, G)$
  under $H$. It has been studied
in \cite{Dijk-Hille-jfa,
gz-bere-rbsd,gz-adv}; 
 see also \cite{Neretin-plan-beredef}
for some general higher
rank cases. Let $X=H/L$ be the 
symmetric space of  $H$
as in \S2. The
 identification of $H$ with 
a subgroup in $G$ induces
an imbedding of $X$ in $D=G/K$
as a totally real totally geodesic
submanifold.

{Let 
$R: \mathcal H_\nu \to C^\infty(X)$
be the restriction
$$RF(x)= \begin{cases} F(x), \ 
x\in X\\
0 \ \mbox{ otherwise}
\end{cases}
$$
Define the weighted restriction  $R_\nu: \mathcal H_\nu \to C^\infty(X)$ by
$$
R_\nu F(x)=(1-\langle x, x\rangle)^{\lambda} RF(x)=
(1-\langle x, x\rangle)^{\lambda}
 F(x), x\in X\subset D, $$
 where $$
\lambda= 
\begin{cases} 
\frac {\nu} 2,  
  & \mathbb F=\br \\
\nu, &  \mathbb F=\bc\\
\nu,  &  \mathbb F=\bh.
\end{cases}
$$
 }
  The functions in 
 $(\mathcal H_\nu, G)$ are holomorphic and thus smooth so $R$ and 
 $R_\nu$ are well-defined. 
 $R_\nu$ is an $H$-intertwining operator
$$
R_\nu: (\mathcal H_\nu, G)
 \to (C^\infty(X), H)$$
with the regular action of  $H  $
on $C^\infty(X)$.

\medskip

We denote by
$\mathcal H_\nu(\bar D)$
 the space of holomorphic
functions on some neighborhoods of $D$; the group
$G$ acts $\mathcal H_\nu(\bar D)$
 via the action $U_\nu$.
In particular elements in $\mathcal H_\nu(\bar D)$
have smooth restrictions to $S$,
 where $S$ is the sphere  in (\ref{eq:1-S-var}), i.e.,
the boundary of $X$. By abuse of notation we denote this restriction also by $R$.
We observe that 
for 
\[ \mu= 
\begin{cases} 
 {\nu},  
  \mbox{ if  \ } \mathbb F=\br \\
2\nu, \mbox{ if  \ }  \mathbb F=\bc, \, \bh,
\end{cases}
\]
$R$ is an $H$-intertrwining operator
\begin{equation}
\label{R-int-rel}
R: (\mathcal H_\nu(\bar D), U_\nu, G)|_H \to 
 (C^\infty (S), \pi_{\mu}, H), \,\,
\end{equation}

If  $\mathbb F= \mathbb R$ this follows immediately from the definition. For the case $\mathbb F=\mathbb C$
we have
\begin{equation*}
\begin{split}
&\quad R(U_\nu (h, h) f)(s)=(U_\nu (h, h) f)(s)\\
&=(cs +d)^{-\nu}  \overline{(cs +d)^{-\nu}} f((as+b)(cs +d)^{-1}, 
\overline{(as+b)(cs +d)^{-1}} )\\
&=|cs +d|^{-2\nu} Rf((as+b)(cs +d)^{-1}), 
\end{split}
\end{equation*}
where $h^{-1}$ is as in  (\ref{eq:H-act}).
The intertwining
relation follows from  (\ref{eq:H-act}). The case of $\mathbb F= \mathbb H$ is  a similar argument using
the observation that $\det (q) =|q|^2$ for any quaternionic
number $q$, with $\det (q)$ being the
determinant of the corresponding $2\times 2$-complex
matrix $q$.

\medskip

Let $R_\nu^\ast$ be the 
adjoint of $R_\nu$ 
considered as a densely defined operator on $L^2(X)=L^2(H/L)$
with the $G$-invariant measure $d\iota(x)$.
Then $R_\nu R_\nu^\ast$ is an integral (generally unbounded) operator 
on  $L^2(X)$
\begin{equation}\label{bere-tr}
R_\nu R_\nu^\ast f(x)=\int_X  B_{\nu}(x, y)       f(y) d\iota(x), \quad B(x, y)=
\left(\frac{(1-\langle x, x\rangle)
(1-\langle x, x\rangle)
}{|(1-\langle x, y\rangle)|^2}\right)^\lambda.
\end{equation}
This operator is called  Berezin transform
\cite{gz-bere-rbsd}, and the representations defined by
the kernel is called the canonical representation \cite{Dijk-Hille-jfa}.
The spectral decomposition of $R_\nu R_\nu^\ast$ is given
in \cite{Dijk-Hille-jfa}.
 As a consequence we
get the branching of  $(U_\nu, G)$  under $H$.
There are finitely many discrete components
and we shall need the first one (compared
in terms of their eigenvalues of the Casimir operator)
of them.

\begin{theo+} Consider
the 
holomorphic representaton $(\mathcal H_\nu(\bar D), U_\nu, G)$. 
Let 
$$\mu=
 \begin{cases} \nu, \,\,       \nu\in (0, \rho_H)=(0, { \frac{n-1}{2})},\quad & \mathbb F=\br  \\
2\nu, \,\,  \nu \in (0, \frac{\rho_H}2)=(0, \frac  n2), \quad  & \mathbb F=\bc \\
2\nu, \,\, \nu\in (1, \frac{\rho_H}{2}) = (1,\frac{2n+1}2), \quad & \mathbb F=\bh.
 \end{cases}
$$
Then the complementary series
$(\pi_{\mu}, H)$ is a discrete summand of the restriction of
in $(\mathcal H_\nu(\bar D), U_\nu, G)$ to $H$. 
Moreover it has  multiplicity one in the restriction. In particular
 the operator
$$R: \mathcal H_\nu(\bar D) \to 
 C^\infty (S)$$ extends a $H$-invariant bounded operator
from   $\mathcal H_\nu$  onto $\pi_{\mu}$.
\end{theo+}

\medskip
Note that the relation between the parameters $\nu$ and $\mu$ is
 determined
by the intertwining property (4.1).

We shall  need concrete constructions
of the intertwining operators.
 The following
result can be derived from \cite[Theorem 5.2]{gz-adv},
where general higher cases are treated.
We present the construction and include 
a proof of the theorem here  in the rank one case.
\begin{theo+}Let  $\nu$  be as in Theorem 4.1.
 Identifying the dual space
of $\pi_{\mu}$ with 
 $\pi_{2\rho- \mu}$  via the standard $G$-invariant paring in $L^2(S)$,
$$
(f, g)=\int_{S} f(x) \overline{g(x)} dx
$$
the dual $J_{\nu} =R^\ast$ of $R: \mathcal H_{\nu}\to
\pi_{\mu}$ is given by
$$
 J_{\nu}: (\pi_{2\rho- \mu}, H)\to  (U_{\nu}, G)|_H, 
\quad 
J_{\nu} 
f(z)=\int_S\frac {f(s)} {h(z, s)^\nu} ds, \, z\in D'
$$
\end{theo+}

{
\medskip

The operator $J_\nu$ is defined as an integral operator with kernel
${h(z, s)^{\nu}}$, the reproducing kernel of $H_{\nu}(D)$. 
We  need its  expansion   \cite{FK} into $K$-types of 
of $H_{\nu}(D)$; 
for the cases of $\mathbb F=\mathbb R, \mathbb C$  it is just the binomial expansion.

Let $\mathcal P(V)$ be the space of polynomials on $V$. It has
\cite{FK-book, FK} an irreducible decomposition 
 under $K$,  known as Hua-Schmid decomposition
$$
\mathcal P(V)
=\begin{cases} \sum_{m=0}^\infty \mathcal P_m(\mathbb C^n)
 & \text{Case} 
\, \mathbb R
\\
\sum_{m_1, m_2=0}^\infty \mathcal P_m(\mathbb C^n)
\otimes \overline{P_m(\mathbb C^n)}
 & \text{Case} 
\, \mathbb C
\\
 \sum_{\m=(m_1, m_2), m_1\ge m_2\ge 0}
 \mathcal P_{\m}
 & \text{Case }
\, \mathbb H
\end{cases}
$$
where  $\mathcal P_m(\mathbb C^n)$ is the space of polynomials of degree 
$m$,
and  for the case $\mathbb H$ the space
  $\mathcal P_{\m}$ is an irreducible
representation of $K$ generated by $z_{11}^{m_1-m_2}\det (z_1)^{m_2}$
for $z=[z_1, \cdots, z_n]^T\in V$
with $z_1\in M_{2, 2}(\mathbb C)$ (of highest
weight $-m_1\gamma_1 -m_2\gamma_2$ with
$\gamma_1, \gamma_2$ being the Harish-Chandra roots).

Denote $(\alpha)_m=\prod_{j=0}^{m-1}(\alpha +j)$
and $(\alpha)_{\m}=\prod_{j=0}^{m_1-1}(\alpha +j)
(\alpha-1 +j)
$ the Pochammer symbol.
Let 
$K_m(z, w)$ and
$K_{\m}(z, w)$
 be the reproducing kernel of the irreducible spaces
$\mathcal P_m$ and 
$\mathcal P_{\m}$ with the Fock space norm.
We have then 
$$K_m(z, w)=\frac{1}{m!}\langle z, w\rangle^m
$$
in Case $\mathbb R$ and
$$
\quad K_{\m}(z, w)=\frac
{d(\m)}
{(2n)_{\m}}
\frac{ \tr_{m_1-m_2}(w^\ast z)}{m_1-m_2+1}
\det(w^\ast z)^{m_2}
$$
in Case $\mathbb H$, the first formula
is elementary
whereas the second
can be obtained
from  \cite[Lemma 3.1]{gz-adv}.
The Faraut-Koranyi expansion
is then the binomial series
$$
h(z,w)^{-\nu}=\sum_{m=0}^\infty (\nu)_m K_m(z, w), \quad K_m(z, w)
=\frac{1}{m!}\langle z, w\rangle^m
$$
for  $\mathbb F= \mathbb R$,
and 
$$
h(z,w)^{-\nu}=\sum_{\m \ge 0}^\infty (\nu)_{\m} 
K_{\m}(z, w),
$$
if $\mathbb F=\mathbb H$. Equivalently this can be rephrased
as 
\begin{equation}\label{fk-pro}
\Vert f\Vert_{\mathcal F}^2 =(\nu)_{\m}  \Vert f\Vert_{\nu}^2.
\end{equation}}
The case $\mathbb F= \mathbb C$ is deduced from the real case.

We can now sketch the proof of the theorem.

\medskip

\begin{proof} 
The formula for $J_\nu$ is obtained 
using the reproducing kernel property. 
Since
 $(\pi_{2\rho- \mu}, H)$ is unitary and irreducible
we need only to check that $J_{\nu}$ maps
the constant function $1$ in $(\pi_{2\rho- \mu}, H)$
into a non-zero element in $(U_{\nu}, G)$, namely
we have to prove
that
\begin{equation}
  \label{eq:j-on-1}
F(z):=J_{\nu}1=\int_S\frac {1} {h(z, s)^{\nu}} ds, \, z\in D  
\end{equation}
is in $ (U_{\nu}, G)$.

$\mathbb F=\br$. Using the power series expansion of $h(z, s)^{-\nu}$ we find
$$
F(z)=\sum_{m}\frac{(\nu)_m}{m!} \int_S \langle z, s\rangle^m ds.
$$
Clearly the integral 
$\int_S \langle z, s\rangle^m ds$
is zero unless  $m=2l$ is even, in which case it
is a constant multiple of the polynomial
$q(z)^l:=(z_1^2 +\cdots +z_n^2)^l$. The constant
can be evaluated at $z=e_1$, namely
\begin{equation*}
  \begin{split}
 \int_S \langle z, s\rangle^{2l} ds
&=q(z)^l 
\int_S s_1^{2l} ds
=C_n B(l+\frac 12, \frac{n-1}2) q(z)^l \\
&=C_n B(l+\frac 12, \frac{n-1}2) q(z)^l
=C_n' \frac{\Gamma(l+\frac 12)}
{\Gamma(l+\frac n2)}
q(z)^l ,
  \end{split}
\end{equation*}
where $C_n$ and $C_n'$ are constants
depending only on $n$. Thus
$$
F(z)=C_n'
\sum_{l}
\frac
{(\nu)_{2l}} 
{(2l)!} 
\frac
{
\Gamma(l+\frac 12)
}
{
\Gamma(l+\frac n2)
}
q(z)^l 
$$
and the membership of $F$  in $(U_\nu, G)$ is
determined by the convergence of square norm
$$
\Vert
F\Vert_\nu^2
=C_n'^2
\sum_{l}
(\frac{(\nu)_{2l}}{(2l)!} 
\frac{\Gamma(l+\frac 12)}
{\Gamma(l+\frac n2)})^2
\Vert q^l \Vert_\nu^2
=C_n'^2
\sum_{l}
(\frac{(\nu)_{2l}}{(2l)!} 
\frac{\Gamma(l+\frac 12)}
{\Gamma(l+\frac n2)}
)^2
\frac{1}{(\nu)_{2l} }\Vert q^l \Vert_{\mathcal F}^2
$$
The square norm $\Vert q^l \Vert_{\mathcal F}^2$
is
$$\Vert q^l \Vert_{\mathcal F}^2
=2^{2l}l! (\frac n2)_l 
=\frac{
\Gamma(\frac 12)
(2l)!
}
{
\Gamma(\frac 12 +l)}(\frac n2)_l.
$$
Simplifying we find
$$
\Vert
F\Vert_\nu^2
=C_n'' \sum_l 
(\frac{(\nu)_{2l}}{(2l)!})
(
\frac{\Gamma(l+\frac 12)}
{\Gamma(l+\frac n2)}
)
\sim \sum_{l\ge 1}\frac{1}{l^{-\nu +\frac n2 +\frac 12}}  
$$
which is convergent if and only 
$\nu<\frac {n-1}2=\rho$.

$\mathbb F=\bc$. This case is
studied in \cite{ehprz-comser}
where a  construction for all discrete
components in the branching of 
$(U_\mu, G)|_{H}$ is obtained. The  function
$F$ above is
$$
F(z, w)=\sum_{m=0}^\infty
\frac{(\nu)_m (\nu)_m}
{(n)_m} \frac 1{m!} 
\langle z, w\rangle ^m, \quad (z, w)\in G/K=B\times B
$$
and 
$$
\Vert F\Vert_{\nu}^2=\sum_{m=0}^\infty
(
\frac{(\nu)_m (\nu)_m}
{(n)_m}
)^2
\frac{\text{dim}(\mathcal P_m)
}{(\nu)_m^2}
$$
where $\text{dim}(\mathcal P_m)\sim m^{n-1}
$ is
the dimension of the polynomial space $P_m(\mathbb C^n)$. 
Here we have used (\ref{fk-pro}) applied to
the tensor product. Each term 
is
$$
\frac
{
(\nu)_m^2
}
{
(n)_m^2
}
{\text{dim}(\mathcal P_m)}
\sim \frac 1{m^{2n-2\nu}}  m^{n-1} 
=
\frac 1{m^{n+1-2\nu}} 
$$
and the series is convergent if an only if
$2\nu <n=\rho_{\fsu(n, 1)}$.

$\mathbb F=\bh$. Similarly $F(z)=\sum_{\m}
(\nu)_{\m}
 \int_S K_{\m}(z, s) ds,
$
and the polynomial
$\int_S K_{\m}(z, s) ds$ 
is an $L$-invariant polynomial in $\mathcal P_{\m}$.
It is proved in \cite[Lemma 3.3]{gz-adv}
that this space is zero unless $\m=(m, m)$, in which
case it is one-dimensional, namely
given by the polynomial $q(z)^m=(\det z_1 +\cdots +\det z_n)^m$.
$$\int_S K_{\m}(z, s) ds
=C_{\m} 
q(z)^m
$$
We evaluate it at $z=[I_2, \cdots, 0]^T$.
The left hand side is
\begin{equation*}
\begin{split}
\int_{S} \det (s_1)^m ds
=\int_{S} |s_1|^{2m} ds
&=C
\int_{\bh} |y|^{2m} (1-|y|^2)^{ 2n-3}
 dy\\
&=C'\int_{0}^1 r^{2m} (1-r^2)^{ 2n-3} r^3 dr
 dy\\
&=C''\frac{(m+1)!}{(2n-1)_m}
\end{split}
\end{equation*}
with a normalization constant $C$ depending only
only on $n$.
The square norm of $F$ 
in
 $(U_{\mu}, G)$ is
$$
\Vert F\Vert_{\nu}^2= C
\sum_{\m=(m, m)} (\nu)_{\m}^2 
\left(\frac{(m+1)!}{(2n-1)_m}\right)^2
\Vert q^m\Vert^2_{\mathcal F}
$$
The Fock space norm $\Vert q^m\Vert^2_{\mathcal F}$
of  $p^m$ can be computed by elementary
differentiation. We have the following
(Caley-Capelli type) identity
$$
q(\partial) q^m=m (m-1+2n) q^{m-1}.
$$
from which it follows
$$
\Vert q^m\Vert_{\mathcal F}^2=m! (2n)_m
$$
The square norm of 
$$
\Vert F\Vert_{\nu}^2
 \sim \sum_{m=1}^\infty \frac 1{m^{2n+2n-1-(2n+2\nu -1)}}
$$
which is convergent precisely when $2\nu < 2n +1
=\rho_{\fsp(n, 1)}$.
\end{proof}
\medskip

\section{Discrete summands in the restriction of 
  $(\pi_\nu, H)$ to the subgroup $H_1$
} 

We prove now our main theorem of this paper.

\begin{theo+}\ \\
A.)  \  Assume that $ \mathbb F= \br$  and $0<\mu <\rho_{H}$. If $\mu+2j <\rho_{H_1}$  $j\ge 0$, then  $(\pi_{\mu+2j} ,H_1) $ is a subrepresentation of the restriction of $(\pi_{\mu},H)$ to $H_1$.\\
B.) \  Assume that $ \mathbb F= \bc$  and $0<\mu  <\rho_{H}$.
If $\mu+2j <\rho_{H_1}$  $j\ge 0$, then  $(\pi_{\mu+2j} ,H_1) $ is a subrepresentation of the restriction of $(\pi_{\mu},H)$ to $H_1$.\\
C.) \ Assume that $ \mathbb F= \bh$  and $2<\mu  <\rho_{H}$. If $\mu+2j <\rho_{H_1}$ $j\ge 0$, then  $(\pi_{\mu+2j} ,H_1) $ is a subrepresentation of the restriction of $(\pi_{\mu},H)$ to $H_1$.
\end{theo+} 

Actually we will prove in the next section
that for  $ \mathbb F= \br$, $0<\mu <\rho_{H}$,
  $\mu+j <\rho_{H_1}$  $j\ge 0$, then  $(\pi_{\mu+j} ,H_1) $ is a subrepresentation of the restriction of $(\pi_{\mu},H)$ to $H_1$. Our proof below
uses explicit differentiations
of $\partial_{z_n^{j'}}$
on $(z_1^2+\cdots +z_{n-1}^2 +z_n^2)^m$ evaluated
at $z_n=0$ which gives trivial result if $j'$ is odd. Presumably
the same proof will also work for odd $j'$ by replacing
the function $1$ below by $x_1$.

\begin{proof} 
Recall that  $\mu$ is defined as a function of $\nu$ by
\begin{equation}
  \label{eq:mu-nu-cor}  
 \mu=
 \begin{cases} \nu, \,\,       \nu\in (0, \rho_H)=(0, { \frac{n-1}{2})},\quad & \mathbb F=\br  \\
2\nu, \,\,  \nu \in (0, \frac{\rho_H}2)=(0, \frac  n2), \quad  & \mathbb F=\bc \\
2\nu, \,\, \nu\in (1, \frac{\rho_H}{2}) = (2,\frac{2n+1}2), \quad & \mathbb F=\bh.
 \end{cases}
\end{equation}
 We consider the operator, recalling 
$\mathcal D^{j} $ in Theorem 3.2,
\begin{equation}
 \label{eq:T-op}
T_j
=(J_{\nu+j'})^\ast 
\mathcal D^{j'} 
J_\nu
=R\,
\mathcal D^{j'} J_\nu,
 \end{equation}
where $j'=2j$ for $ \mathbb F=\mathbb R$
and $j' =j$ for $ \mathbb F=\mathbb C, \mathbb H$.
(We choose $j'=2j$ 
for $ \mathbb F=\mathbb R$ to ensure
that
$(T_j 1, 1)$ below is non-zero.)
 Note that it should be clear
from the formula the spaces 
on which the $J$-operators are defined.
The operator $T_j$ is bounded
and $H_1$-intertwining
$$
T_j=: (\mathcal C_{2\rho_{H}-\nu},  \pi_{2\rho_{H}-\nu}, H)
\stackrel{J_\nu}{\to}
 (\mathcal H_{\nu}, U_\nu, G)
\stackrel{\mathcal D^{j'}}{\to} 
(\mathcal H_{\nu+j'}, U_{\nu+j'}, G_1)
\stackrel{R}{\to} (\mathcal C_{\mu(j')}, \pi_{\mu(j')}, H_1)
$$
by Theorems 3.2 and 4.2, where $\mu(j)$  is obtained 
from $\nu+j'$ via the above correspondence $\mu\leftrightarrow \nu$.
Composed  with  the equivalence 
$\Lambda_{2\rho_{H}-\nu}$ of
$(\mathcal C_{\nu},  \pi_{\nu}, H)$
and $(\mathcal C_{2\rho_{H}-\nu},  \pi_{2\rho_{H}-\nu}, H)$
we get a bounded $H_1$-intertwining operator
$$
T_j\Lambda_{\nu}:
(\mathcal C_{\nu},  \pi_{\nu}, H)\to 
(\mathcal C_{\mu(j)}, \pi_{\mu(j)}, H_1)
$$
which preserves the Hilbert space norms on $H$ and $H_1$.
Thus we need only
to check that $T_j1\ne 0$. We compute
the pairing
$(T_j 1, 1) $ in  Theorem 4.2.
$$
(T_j 1, 1)
=(\mathcal D^{j'}
J_\nu 1,
 J_{2\rho_{H_1}-\mu(j)}
1
)_{\mathcal H_{\nu+j'}}.
$$
The functions $J_\nu 1$
and $J_
{
2\rho_{H_1}-\mu(j)
} 1$ are found in 
  (\ref{eq:j-on-1} ) in the proof of Theorem 4.2
(for the bigger group pair $(G, H)$).
 We shall find
 $\mathcal D^k
J_\nu 1$ following the computations there.

$\mathbb F=\mathbb R$. In this case
$\mathcal D^{2j} 
J_\nu 1$
is, apart from a positive normalization constant,
$$
\sum_{l\ge j}
\frac
{(\nu)_{2l}} 
{(2l)!} 
\frac
{
\Gamma(l+\frac 12)
}
{
\Gamma(l+\frac n2)
}
(2j)!\binom l{2j} 
q(z', 0)^{l-j}
$$
for 
$$
(\partial_n^j q^l)(z',0) =
(2j)!\binom l{2j} 
q(z', 0)^{l-j},
$$
whereas 
 $J_{2\rho_{H_1}-\mu(j)} 1$ is
$$
F(z')=C
\sum_{k}
\frac
{(2\rho_{H_1}-\mu(j))_{2k}} 
{(2k)!} 
\frac
{
\Gamma(k+\frac 12)
}
{
\Gamma(k+\frac n2)
}
q(z')^k.
$$
The inner product $
(\mathcal D^{2j}
J_\nu 1,
 J_{\mu} 1
)_{\mathcal H_{2\rho_{H_1}-\mu(j)}}
$, is then,
up to a positive constant,
$$
\sum_{l-j=k=0}^\infty
\frac
{(\nu)_{2l}} 
{(2l)!} 
\frac
{
\Gamma(l+\frac 12)
}
{
\Gamma(l+\frac n2)
}
(2j)!\binom l{2j} 
\frac
{(2\rho_{H_1}-\mu(j))_{2k}} 
{(2k)!} 
\frac
{
\Gamma(k+\frac 12)
}
{
\Gamma(k+\frac n2)
}\langle q(z')^k, q(z')^k
\rangle_{\mathcal H_{\mu(j)}},
$$
which is a (convergent) sum
of positive terms and thus is non-zero.

$\mathbb F=\mathbb C,\mathbb H$. The computations
are similar. We need only to note that for the $\mathbb H$-case,
$$
\partial_Q^j q^m (z', 0)=\binom{m}{j} (j+1)! j!\, q(z', 0)^{m-j}
$$
with the coefficients
being also positive.
\end{proof}

We note that  our approach
above can also be used to treat
the complementary series
for $H=F_{4(-20)}$ 
\cite{Johnson-2}
under the symmetric
subgroup $H_1=Spin(8, 1)$ by considering
the Hermitian Lie groups $G\supset G_1$
with $G=E_{6(-25)}$,
$G_1=Spin(8, 2)$.  
However
no more extra discrete component
will be discovered except 
the one 
obtained using the restriction 
method in \cite{gz-brch-rk1},
namely
the zero degree of differentiation 
in the above method. This
is due to the rather short interval
of complementary series of 
 $H=F_{4(-20)}$. More precisely
the eligible range of the parameter $\nu$ of
holomorphic representations
$(U_\nu, G)$
 for which  discrete
components under $H$ do appear is the interval $(3, \frac{11}2)$.
Here the lower bound $3$ corresponds
to the first discrete Wallach reducible
point for $G$,  and the upper bound
corresponds to $\frac{\rho_{H}}2$
 which is the critical point where the
spectral symbol of the Berezin 
transform $RR^\ast$ in (\ref{bere-tr})
has a pole \cite[(6.1)]{gz-brch-rk1}. The same
range 
of $\nu$ for $(U_\nu, G_1)$
for the pair $(G_1, H_1)=(Spin(8, 2), Spin(8, 1))$
is $(3, \frac 72)$. 
(The $\rho_{H}$ and $\rho_{H_1}$
are $11$ and $7$ with the same normalization of
split Cartan in $\fh_1\subset \fh$, the half values
$\frac{11}2$ and $\frac 72
$ in the intervals are due to the
double rank of $G$ and $G_1$
of that of $H$ and $H_1$.)
The symmetric domain $D_1=G_1/K_1$ is in $V_1=\mathbb O_{\mathbb C}
$
as a subdomain of $D=G/K$  in $\mathbb O_{\mathbb C}^2=
V_1 \oplus V_2$.
A normal differentiation by any $K_1$ invariant
polynomial on $V_2$ on the subdomain $D_1=G_1/K_1$
will result
in intertwining operator
from $(U_\nu, G)$ into
$(U_{\nu'}, G_1)$
with $\nu' \ge \nu +1$. But for $\nu >3$
we have $\nu' \ge 4 >\frac 72$
which is outside the eligible range 
$(3, \frac 72)$.
Presumably
this is the only possible
discrete complement.

\section{Explicit realization
of discrete components via
intertwining differential
operators
}

For the real case a direct proof
of the appearance of the discrete
components in the branching
can be obtained using the
explicit formula for intertwining
differential operators \cite{Juhl-book2}
for
spherical representations
of $H$ and $H_1$. The operators
 have an easier expression in the non-compact
realization. In this realization
our proof is a simple application of
the Cauchy-Schwarz inequality, and it also yields
an explicit construction of the discrete components.

The complementary series 
$(\mathcal C_\mu, \pi_{\mu}, H)$, 
$H=SO_0(n, 1)$,
can be realized as distribution space on $\fn=\mathbb R^{n-1}$; see e.g. \cite{Speh-Venk-2}. Initially it  acts on 
a subspace 
$C^\infty_\mu(\mathbb R^{n-1})$
of $C^\infty(\mathbb R^{n-1}
)$ obtained from the restriction
to $\mathbb R^{n-1}=N^- 
$ of smooth functions on $G$;
alternatively it can be obtained by a
 a Cayley transfrom from $S$ to $\mathbb R^{n-1}$
and the $H$-action on the space 
$C^\infty(S)$. The action of $P^-=MAN^-$ on 
the space $C^\infty_\mu(\mathbb R^{n-1})$ 
is
$$
\pi_{\mu}(m)f(x)= f(m^{-1}x), \,\,
\pi_{\mu}(a)f(x)= a^{-\mu}f(a^{-1}x), \,\,
\pi_{\mu}(y)
f(x)= f(x-y), 
$$
for 
$$
\, m\in SO_0(n-1),\,
\, a\in\mathbb R^+=A,  \,\,
 y\in \mathbb R^{n-1}=
N^-.
$$
The Weyl group element $w$ acts as
an inversion
$$
\pi_{\mu}(w)f(x)=
f(-\frac{2x}{|x|^2})
(\frac{2}{|x|^2})^{\mu}
$$ 
The Knapp-Stein intertwining operator
is
$$
f\to \int_{\mathbb R^{n-1}} f(y) |x-y|^{-2(n-1-\mu)} dy 
$$
as a meromorphic continuation of integral operators.
Let $\Delta_{n-1}=\partial_1^2 + \cdots +\partial_{n-1}^2$ be
the Laplace operator on $\mathbb R^{n-1}$
and $\mathcal F$ the Fourier transform. 
Then $(\mathcal C_\mu, \pi_{\mu}, H)$ is the completion of 
$C_0^\infty(\mathbb R^{n-1})$
with the norm
\begin
{equation}
\label{norm-cpl}
\Vert f\Vert_{\mu}^2=\Vert
 \Delta^{\frac{n-1-2\mu}2} f\Vert_{L^2}^2 =\int_{\mathbb R^{n-1}}
|\mathcal F f(\xi)|^2 |\xi|^{n-1-2\mu}d\xi
\end
{equation}
for $0<\mu <\frac{n-1}2=\rho$. See e.g. \cite[Theorem 2.1]{Vershik-Graev},
our parameter $\mu$ being their $n-1-\lambda$.

 The explicit form
of $\pi_{\mu}$ is well-known and will not be used
here. 
 Correspondingly
the representation $(\mathcal C_\mu, \pi_{\mu}, H_1)$
is realized on $\mathbb R^{n-2}$, considered
as the subspace $\{(x, 0); x\in\mathbb R^{n-2}\}$
of $\mathbb R^{n-1}$.

We recall the result of Juhl
\cite[Ch. V]{Juhl-book2} 
on intertwining
differential operators, with our parametrization
of the spherical representations.

\begin{theo+} Let $n>2$  
and $\mu\in \mathbb C$. There exists an non-zero intertwining $H_1$ operator 
$$
\mathcal D_{\mu, j}: 
(C^\infty_\mu(\mathbb R^{n-1}), \pi_{\mu}, H)
\to (C^\infty_{\mu+j}(\mathbb R^{n-2}), \pi_{\mu+j}, H_1)
$$
of the form
\begin{equation}
\label{D-op}
\mathcal D_{\mu, j} f(y,0)
=\left
(
\sum_{2k\le j} C_{\mu, k}\partial_{n-1}^{j-2k} \Delta_{n-2}^{k} f
\right)
(y, 0)
\end{equation}
where $C_{\mu, k}$
are constants depending on $\mu, m$ normalized
by $C_{\mu, 0}=1$.
\end{theo+}

We give now an independent  proof of Theorem 5.1 for the real case.

\begin{theo+} Let  $n>2$,  $0<\mu<\rho_H=\frac{n-1}2$ and 
$\mu +j <\rho_{H_1}=\frac{n-2}2$.
Then the operator $\mathcal D_{\mu, j} $
is a non-zero $H_1$-invariant
bounded operator from  $(\mathcal C_\mu, H)$
onto  $(\mathcal C_{\mu+j}, H_1)$. 
Thus  the complemetnary series $(\mathcal C_{\mu+j}, H_1)$
appears as a discrete component in  the decompostion of 
$(\mathcal C_\mu, H)$
 under $H_1$.
\end{theo+} 
\begin{proof} 
Note that the second claim follows from the first
one by a standard argument. Indeed
the adjoint $ \mathcal D_{\mu, j}^\ast$
of 
 $\mathcal D_{\mu, j} $
is then a non-zero bounded $H_1$-invariant
operator from
 $(\mathcal C_{\mu+j}, H_1)$
to $(\mathcal C_\mu, H)$. Its polar decomposition
defines an $H_1$-invariant
partial isometry. But since 
 $(\mathcal C_{\mu+j}, H_1)$
is irreducible it is an isometry.

We denote $\mathcal F^\flat$
the Fourier transform on $\mathbb R^{n-2}$.
Let $f\in C_0^\infty(\mathbb R^{n-1})$
and
$g:=D_{\mu, j} 
f$. We shall prove that
$$
\Vert g
\Vert_{\mu +j}^2
:=\int_{\mathbb R^{n-2}}
|\eta|^{n-2-2(\mu+j)}|\mathcal F^\flat
g(\eta)|^2 d\eta
\le C \Vert f\Vert_{\mu }^2.
$$
By the Fourier inversion formula
we have, writing the dual variable $\xi\in \mathbb R^{n-1}$ 
as $\xi=(\eta, \xi_{n-1})\in \mathbb R^{n-2}\times \mathbb R$, 
$$
f(y, x_{n-1})=
C_n\int_{\mathbb R^{n-1}}
\mathcal F f(\eta, \xi_{n-1}) e^{i(y, \eta) +i x_{n-1} \xi_{n-1}} d\eta d\xi_{n-1},
$$
where $C_n$ is the normalization constant for the Fourier inversion.
Performing the differentiation 
$\mathcal D_{\mu, j} $ evaluated at $(y, 0)$ we get
\begin{equation*}
\begin{split}
g(y)&=\mathcal D_{\mu, j} f(y, 0)
=C_n\int_{\mathbb R^{n-1}}
\mathcal Ff(\eta, \xi_{n-1}) e^{i(y, \eta) } E(\eta, \xi_{n-1})d\eta d\xi_{n-1}\\
&=
C_n\int_{\mathbb R^{n-2}}e^{i(y, \eta) } \left(
\int_{\mathbb R} 
\mathcal Ff(\eta, \xi_{n-1}) 
E(\eta, \xi_{n-1})
d\xi_{n-1}   \right) 
d\eta   
\end{split}
\end{equation*}
where $E(\xi)=E(\eta, \xi_{n-1})$  is the Fourier
multiplier corresponding to (\ref{D-op}),
$$
E(\xi)=E(\eta, \xi_{n-1})=\sum_{2k\le j}
C_{\mu, k}(-|\eta|^2)^{k} (-\xi_{n-1}^2)^{j-2k} 
.$$
 The Fourier transform $\mathcal F^\flat g$ of $g$
is
$$
\mathcal F^\flat g(\eta)= 
C
\int_{\mathbb R}  
 \mathcal Ff(\eta, \xi_{n-1}) 
E(\eta, \xi_{n-1})
d\xi_{n-1}
=C
\int_{\mathbb R}  
 \mathcal Ff(\xi) 
E(\xi)
d\xi_{n-1}
$$
where we have written $\xi=(\eta, \xi_{n-1})$.
We apply the Cauchy-Schwarz inequality, obtaining
$$
|\mathcal F^\flat g(\eta)|^2 
\le 
\int_{\mathbb R} 
| \mathcal F f(\xi) |^2 |\xi|^{n-1-2\mu} 
d\xi_{n-1}
\int_{\mathbb R} 
|E(\xi)|^2 
 \frac 1{|\xi|^{n-1-2\mu} }
d\xi_{n-1}.
$$
Now the polynomial $E(\xi)$ is a sum of terms of the form $\xi_{n-1}^{j-2k} |\eta|^{2k}$, $0\le 2k\le j$,  and we apply the triangle inequality to 
the second integration,
$$
\int_{\mathbb R} 
|E(\xi)|^2 
 \frac 1{|\xi|^{n-1-2\mu} }
d\xi_{n-1}
\le C
\sum_{k=0}
\int_{\mathbb R} 
(|\eta|^{2k} \xi_{n-1}^{j-2k})^2
 \frac 1{|\xi|^{n-1-2\mu} } d\xi_{n-1}.
$$
The each term can be evaluated by a change of variables,
$\xi_{n-1}=|\eta| u$, 
$$
\int_{\mathbb R} 
(|\eta|^{2k} \xi_{n-1}^{j-2k})^2
 \frac 1{|\xi|^{n-1-2\mu} }
d\xi_{n-1} =
\frac 1{|\eta|^{n-2-2\mu-2j} }
\int_{\mathbb R} \frac {u^{4k}}
{(1 +u^2)^{n-1-2\mu}} du
=C\frac 1{|\eta|^{n-2-2\mu-2j} }.
$$
Here we have used
the assumption that $j\le j_0$ so that $n-1-2\mu-4k\ge
n-1-2\mu-2j >1$
and the integration 
$$
C=\int_{\mathbb R} \frac {u^{4k}}
{(1 +u^2)^{n-1-2\mu}} du <\infty
$$
is convergent. That is
$$
\int_{\mathbb R} 
|E(\xi)|^2 
 \frac 1{|\xi|^{n-1-2\mu} }
d\xi_{n-1} \le C\frac 1{|\eta|^{n-2-2\mu-2j} }
$$
and
$$
|\mathcal F^\flat g(\eta)|^2 
\le C
\frac 1{|\eta|^{n-2-2\mu-2j} }
\int_{\mathbb R} 
| \mathcal Ff(\xi) |^2 |\xi|^{n-1-2\mu} 
d\xi_{n-1}.
$$
Integrating over $\eta\in \mathbb R^{n-2}$
we find
$$
\Vert g\Vert^2_{\mu+j}=\int_{\mathbb R^{n-2}}
|\mathcal F^\flat g(\eta)|^2 
|\eta|^{n-2-2\mu-2j} 
\le C
\int_{\mathbb R^{n-2}} \int_{\mathbb R} 
|\mathcal Ff(\xi) |^2 |\xi|^{n-1-2\mu} d\xi
=C \Vert f\Vert_{\mu}^2
$$
completing the proof.
\end{proof}

For $j=0$, i.e.  $\mathcal D_{\mu, 0} $ being the restriction,
 the above result is proved in \cite{Speh-Venk-2} by considering
the dual map of the restriction. Note the above
proof is somewhat simpler.

Note the standard $L^2(\mathbb R^{n-1})$-pairing 
 $f\otimes g\to (f, g)$ 
is invariant under  $\pi_{\mu}\otimes \pi_{2\rho-\mu}$
and we may rephrase the above theorem as
a realization of the discrete components.
\begin{coro+} Let $\mu$ and $j$ be as above.
The representation $(\pi_{n-2-(\mu +j)}, H_1)$
appears in the branching of
$(\pi_{n-1-\mu}, H)$ under $H_1$ and
an explicit realization
of the subrepresentation is the distribution
space
$$
\{F; 
F(y, x_{n-1})=
\sum_{2k\le j} C_{\mu, k}
\partial_{n-1}^{j-2k} \Delta_{n-2}^{k} 
(f(y) \delta(x_{n-1})), \quad f\in (\pi_{n-2-(\mu +j)}, H_1)
\},
$$
where $\delta$ is the Dirac's delta function.
\end{coro+} 

We may also use the version of
$\mathcal D_{\mu, j}$ in the compact
picture and to prove Theorem 6.2 directly
using the method in \cite{gz-brch-rk1}.
This involves rather detailed computations
of normal differentiation of spherical harmonics.
The above proof using Fourier transform
is apparently easier.

Finally we study the limit case $\mu=0$ using the operator
$\mathcal D_{\mu, j} $. Assume $n >2$ is an odd integer.
The Harish-Chandra module
of $K$-finite elements
in 
$C^\infty(S)$ has the constant functions $\mathbb C$
 as its submodule 
and the quotient space is unitarizable,
 with the
unitary norm obtained from that of 
$\mathcal C_\mu$ by continuation (or residue depending
on the normalization); see \cite{Johnson-Wallach}.
We denote the quotient space by $\mathcal W_0$.
Now the  norm
\begin{equation}
\label{norm-cpl-quo}
\Vert f\Vert_{0}^2=\Vert
 \Delta^{\frac{n-1}2} f\Vert_{L^2}^2 =\int_{\mathbb R^{n-1}}
|\mathcal F f(\xi)|^2 |\xi|^{n-1}d\xi
\end
{equation} 
is the continuation of 
(\ref{norm-cpl}) at $\mu=0$ and it vanishes
on the constant function. Thus it is the unitary norm
on $\mathcal W_0$. See also \cite[Theorem 2.2]{Vershik-Graev}.
Using the same proof
as above
we have the following result
 \begin{theo+} Let  $n>2$.
Then the operator $\mathcal D_{0, j} $
is a bounded from  $\mathcal W_0$
onto  $(\mathcal C_{j}, H_1)$ for $1\le j<\frac{n-2}2$. Thus
$(\mathcal C_{j}, H_1)$, $1\le j<\frac{n-2}2$,
appear as a discrete component of
$(\mathcal W_0, H)_{H_1}$.
\end{theo+}

\def\cprime{$'$} \newcommand{\noopsort}[1]{} \newcommand{\printfirst}[2]{#1}
  \newcommand{\singleletter}[1]{#1} \newcommand{\switchargs}[2]{#2#1}
  \def\cprime{$'$} \def\cprime{$'$} \def\cprime{$'$}
\providecommand{\bysame}{\leavevmode\hbox to3em{\hrulefill}\thinspace}
\providecommand{\MR}{\relax\ifhmode\unskip\space\fi MR }
\providecommand{\MRhref}[2]{%
  \href{http://www.ams.org/mathscinet-getitem?mr=#1}{#2}
}
\providecommand{\href}[2]{#2}


\end{document}